\documentclass[11pt,
]{article}
\usepackage{amssymb,amstext,amsmath,amsthm,latexsym,mathrsfs,mathbbol}

\newcommand{\nek}{\newcommand}

\nek{\vyk}[1]{}

\nek{\itsep}{\itemsep=0.4ex plus 0.15ex minus 0.15ex}
\nek{\tenu}[1]{
\def\theenumi{#1}\itsep}

\theoremstyle{plain}

\newtheorem{theorem}             {Theorem}
\newtheorem{corollary}  [theorem]{Corollary}
\newtheorem{lemma}      [theorem]{Lemma}
\newtheorem{proposition}[theorem]{Proposition}

\theoremstyle{definition}

\newtheorem{defi}       [theorem]{Definition}
\newtheorem*{prof}{Proof}

\nek{\TF}{\sl}
\nek{\bcor}{\begin{corollary}}
\nek{\ecor}{\end{corollary}}
\nek{\bdf} {\begin{defi}}
\nek{\edf} {\qed\end{defi}} 
\nek{\eDf} {\end{defi}} 
\nek{\ble} {\begin{lemma}}
\nek{\ele} {\end{lemma}}
\nek{\bpf} {\begin{prof}}
\nek{\epf} {\qed\end{prof}}
\nek{\bpro}{\begin{proposition}} 
\nek{\epro}{\end{proposition}} 
\nek{\bte} {\begin{theorem}\TF\ }
\nek{\ete} {\end{theorem}}

\nek{\qeD} {\qed}
\nek{\qeG} [1]{\hfill{\hbox{$\mtho\square$~({\sl#1\hspace{0.2ex}})}}}
\nek{\qedD}{\hfill{$\msur\dashv\msur$}}
\nek{\ben}{\begin{enumerate}\itsep}
\nek{\een}{\end{enumerate}}
\nek{\bit}{\begin{itemize}\itsep}
\nek{\eit}{\end{itemize}}
\nek{\bay}{\begin{array}}
\nek{\eay}{\end{array}}
\nek{\Da}{\Delta}
\nek{\sg} {\sigma}
\nek{\om} {\omega} 
\nek{\Om} {\Omega} 
\nek{\bbb}{\hspace{0.5pt}} 
\nek{\dvoj}[1]{{\bbb{\mathbb #1}\bbb}}
\nek{\dH}{{\dvoj H}}
\nek{\dQ}{{\dvoj Q}}
\nek{\dR}{{\dvoj R}}
\nek{\dN}{{\dvoj N}}
\nek{\dT}{{\dvoj T}}
\nek{\dZ}{{\dvoj Z}}
\nek{\fsg} {\sg}
\nek{\skrsp}{\hspace{0.5pt}} 
\nek{\scri}[1]{\skrsp\mathscr #1\skrsp} 
\nek{\cP}{{\scri P}}
\nek{\cZ}{{\scri Z}}
\nek{\sneq}{\subsetneqq}
\nek{\sq}  {\subseteq}
\nek{\pu}  {\varnothing}
\nek{\imp} {\Longrightarrow} 
\nek{\res} {\restriction} 
\nek{\ti}  {\times}
\nek{\dm}  {$$}
\nek{\kaz} {\forall\,}
\nek{\sus} {\exists\,}
\nek{\nin} {\not\in}
\nek{\bez} {\smallsetminus}
\nek{\onto}{\stackrel{\rm onto}\lra}
\nek{\ang} [1]{\langle #1\rangle}
\nek{\ans} [1]{\{\hspace{0.2mm}#1\hspace{0.2mm}\}}
\nek{\sis}[2] {{\ans{#1}}_{#2}} 
\nek{\dd} [1] {$\rsur\mtho#1\qsur$-}
\nek{\imar}[1]{\marginpar%
{\vspace{-1ex}\flushleft\footnotesize\it#1}}
\nek{\itla}{\item\label}
\nek{\mtho}{\mathsurround=0mm}
\nek{\msur}{\hspace*{-1\mathsurround}}
\nek{\qsur}{\hspace{0.2\mathsurround}}
\nek{\rsur}{\hspace{0.4\mathsurround}}
\nek{\noi}{\noindent}
\nek{\vom}{\vspace{1mm}}
\nek{\vtm}{\vspace{2mm}}
\nek{\hm}{homomorphism}
\nek{\ddi}{\dd\cI} 
\nek{\ddz}{\dd\cZ} 

\nek{\sd}{\mathbin{\Da}}

\nek{\gd}{{\bf G}_\da}
\nek{\fs}{{\bf F}_\fsg}

\nek{\opl}{\oplus}
\nek{\omi}{\ominus}

\nek{\zz} {\linebreak[0]} 
\nek{\ens} [2] {\ans{{#1\hspace{0.5ex}{:}}\zz\hspace{0.5ex}#2}}

\nek{\aenu}{\tenu{{\rm(\arabic{enumi})}}}
\nek{\renu}{\tenu{{\rm(\roman{enumi})}}}

\nek{\ve} {\varepsilon}

\nek{\qeDD} [1] 
{\hfill\hbox{\qed~({\small #1\/}\hspace{0.1ex})}}

\nek{\epF} [1] {\qeDD{#1}\end{prof}} 

\nek{\kon} {FIN-}

\nek{\bbo} {\mathbb0}

\nek{\raz} [2] {\Da_{#1\hspace*{0.15ex},\hspace*{0.3ex}#2}}

\makeatletter
\if@twoside%
\else\fi%
\makeatother

\begin{document}                                 

\title{Reasonable non--Radon--Nikodym idealss}
\author{Vladimir Kanovei\thanks
{IITP, Moscow.
Partially supported by RFFI grants 06-01-00608 and 07-01-00445.}
\and
Vassily Lyubetsky\thanks
{IITP, Moscow. 
Partially supported by RFFI grant 07-01-00445. }
}
\date{\today}
\maketitle

\begin{abstract}
For any abelian Polish \dd\sg compact group $\dH$ 
there exist an $\fs$ ideal $\cZ\sq\cP(\dN)$ and a Borel 
\dd\cZ approximate homomorphism $f:\dH\to\dH^\dN$ which  
is not \dd\cZ approximable by a continuous true 
homomorphism $g:\dH\to\dH^\dN$.              
\end{abstract}

\subsubsection*{Introduction}

Let $G,H$ be abelian Polish groups, 
and $\cZ$ be an ideal over a countable set $A$. 
We consider $H^A$ as a product group. 
For $s,t\in H^A$ put 
\dm
\raz{s}{t}=\ens{a\in A}{s(a)\ne t(a)}\,.
\dm 
Suppose that $\cZ$ is an ideal over $A$. 
A map $f:G\to H^A$ is a {\it\dd\cZ approximate \hm\/} iff 
$\raz{f(x)+f(y)}{f(x+y)}\in\cZ$ for all $x,y\in G$.
Thus it is required that the set of all $a\in A$ such that 
$f_a(x)+f_a(y)\ne f_a(x+y)$ belongs to $\cZ$. 
Here $f_a:G\to H$ is the \dd ath co-ordinate map of the map 
$f:G\to H^A$. 

And $\cZ$ is a {\it Radon--Nikodym\/} ideal 
(for this pair of groups) 
iff for any measurable \dd\cZ approximate \hm\ $f:G\to H^\dN$ 
there is a continuous exact \hm\ $g:G\to H^\dN$ which  
\dd\cZ approximates $f$ in the sense that
$\raz{f(x)}{g(x)}\in\cZ$ for all $x\in G$.
Here the measurability condition 
can be understood as Baire measurability, 
or, if $G$ is equipped with a \dd\sg additive Borel measure, as 
measurability with respect to that measure. 

The idea of this (somewhat loose) concept is quite clear: the 
Radon--Nikodym ideals are those which allow us to approximate 
non-exact \hm s by true ones. 
This type of problems appears in different domains of mathematics. 
Closer to the context of this note, 
Velickovic \cite{v} proved that any Baire-measurable 
\kon approximate Boolean-algebra automorphism $f$ of $\cP(\dN)$ 
(so that the symmetric differences between $f(x)\cup f(y)$ 
and $f(x\cup y)$ and between $f(\dN\bez x)$ and $\dN\bez f(x)$ 
are finite for all $x,y\sq\dN$) is \kon approximable by a true 
automorphism $g$ induced by a bijection betveen two cofinite 
subsets of $\dN$. 
Kanovei and Reeken proved that any Baire measurable 
\dd\dQ approximate homomorphism 
$f:\dR\to\dR$ is \dd\dQ approximable by a homomorphism of the 
form $f(x)=cx$, $c$ being a real constant.
See also some results in \cite{f,kr1,kr2}.

The term ``Radon--Nikodym ideal'' was introduced by 
Farah~\cite{f,aq} in the context of Baire measurable 
Boolean algebra \hm s of $\cP(\dN)$.
Many known Borel ideals were demonstrated to be Radon--Nikodym, 
see \cite{f,aq,kr1,kr2}.
Suitable counterexamples, again in the context of Boolean 
algebra \hm s, were defined by Farah on the base of so called 
pathological submeasures. 
A different and, perhaps, more transparent counterexample, 
related to \hm s $\dT\to\dT^\dN$ 
(where $\dT=\dR/\dN$), 
is defined in \cite{kr2} as a modification of an ideal 
introduced in \cite{so}.  
The next theorem generalizes this result. 

\bte
\label m
Suppose that\/ $\dH$ is an uncountable abelian Polish group. 
Then there is an analytic ideal\/ $\cZ$ over\/ $\dN$ that is 
not a Radon -- Nikodym ideal for maps\/ $\dH\to\dH^\dN$ in the 
sense that there is a Borel and\/ \dd\cZ approximate \hm\/ 
$f:\dH\to\dH^\dN$ not\/ \dd\cZ approximable by a 
continuous \hm\/ $g:\dH\to\dH^\dN$. 
If moreover\/ $\dH$ is\/ \dd\sg compact then\/ $\cZ$ can be 
chosen to be\/ $\fs$.
\ete

Note that the theorem will not become stronger if we require 
$g$ to be only Baire-measurable, or just measurable with respect 
to a certain Borel measure on $\dH$ --- because by the Pettis 
theorem any such a measurable group \hm\ must be continuous.

The remainder of the note contains the proof of Theorem~\ref m.
It would be interesting to prove the theorem for non-abelian 
Polish groups. 
(The assumption that $\dH$ is abelian is used in the proof 
of Lemma~\ref{l1}.) 
And it will be interesting to find non--Radon--Nikodym ideals 
for \hm s $G\to H^\dN$ in the case when the Polish groups 
$G$ and $H$ are not necessarily equal.

\subsubsection{Countable subgroup}
\label{pr}

Let us fix a group $\dH$ as in the theorem, that is, an 
uncountable abelian Polish group.
By $\bbo$ we denote the neutral element, 
by $\opl$ the group operation, 
by $d$ a compatible complete separable distance 
(and we do not assume it to be invariant).
The first step is to choose a certain countable subgroup 
$D\sq\dH$ of ``rational elements''.

It is quite clear that there exists 
a countable dense subgroup $D\sq\dH$ 
satisfying the following requirement of 
elementary equivalence type. 

\bit
\item[$(\ast)$]
Suppose that $n\ge 1$, $c_1,\dots,c_n\in D$, 
$\ve$ is a positive rational, $U_i=\ens{x\in\dH}{d(x,c_i)\le\ve}$, 
and $P(x_1,\dots,x_n)$ is a finite system of linear equations 
with integer coefficients, unknowns $x_1,\dots,x_n$, and constants 
in $D$, of the form:
\dm
b_1x_1\opl\dots\opl b_nx_n=r\,,\quad \text{where}\quad
b_i\in\dZ\;\text{ and }\;r\in D\,.
\dm
Suppose also that this system $P$ has a solution 
$\ang{x_1,\dots,x_n}$ in $\dH$ such that $x_i\in U_i$ for all $i$. 
Then $P$ has a solution in $D$ as well.
(That is, all $x_i$ belong to $D\cap U_i$.)
\eit

Let us fix such a subgroup $D$.

\subsubsection{The index set}
\label{is}

Let 
{\it rational ball\/} mean any subset of $\dH$ of the form
$\ens{x\in\dH}{d(c,x)<\ve}$, where $c\in D$ (the center), 
and $\ve$ is a positive rational number.

\bdf
\label{3a} 
Let $A$, the index set, consist of all objects $a$ of the 
following kind. 
Each $a\in A$ consists of:
\bit
\itsep
\item[$-$] 
an open non-empty set $U^a\sneq\dH$, 

\item[$-$] 
a partition $U^a=U^a_1\cup\dots\cup U^a_n$ of $U^a$ onto 
a finite number $n=n^a$ of pairwise disjoint non-empty 
rational balls $U^a_i\sq\dH$, and 

\item[$-$] 
a set of points $r^a_i\in U^a_i\cap D$ 
such that, for all $i,\,j=1,2,\dots,n$:
\ben
\aenu
\itla{31} 
{\it either\/} $r^a_i\opl r^a_j$ is $r^a_{k}$ for some $k,$ 
and $(U^a_i\opl U^a_j)\cap U^a\sq U^a_k$,

\itla{32}
{\it or\/} $(U^a_i\opl U^a_j)\cap U^a=\pu$.\qed 
\vyk{\item[$-$] 
{\it or\/} $i=j$ and $s_i=0$.  }
\een
\eit
\eDf

Under the conditions of Definition~\ref{3a}, 
if $\bbo\in U^a_i$ then $s_i=\bbo$: for take $j=i$. 


\ble
\label{ai}
$A$ is an infinite (countable) set.
\ele
\bpf
For any $\ve>0$ there is $a\in A$ such that 
$U^a$ a set of diameter $\le\ve$: 
just take $n^a=1$, $r^a_1=\bbo$, and let 
$U^a=U^a_1$ be the \dd{\frac\ve2}nbhd of $\bbo$ in $\dH$.
\epf

The next lemma will be used below.

\ble
\label{l*}
If\/ $y_1,\dots,y_n\in\dH$ are pairwise distinct then there 
exists\/ $a\in A$ such that\/ $n^a=n$ and\/ 
$y_i\in U^a_i$ for all\/ $i=1,\dots,n$.
\ele
\bpf
As the operation is continuous, we can pick pairwise disjoint 
rational balls $B_1,\dots,B_n$ such that $y_i\in B_i$ for all 
$i$ and the following holds: 
If $1\le i,j\le n$ then either there exists $k$ such that 
$(B_i\opl B_j)\cap B\sq B_k$, where $B=B_1\cup\dots\cup B_n$, 
or just $(B_i\opl B_j)\cap B=\pu$. 
Put $U^a_i=B_i$. 

To obtain a system of points $r^a_i$ required, 
\vyk{
Consider the finite set 
$$
X=\ens{y_i}{1\le i\le n} 
\cup\ens{y_i\opl y_j}{1\le i,j\le n}\,.
$$
Let $\ve$ be the smallest positive real of the form 
$d(x,y)$, where $x,y$ are two distinct elements in $X$.
Put $n^a=n$ and let each $U^a_i$ ($1\le i\le n$) 
be a rational ball in $\dH$ of diameter $\frac\ve4$ 
whose center $c_i\in D$ is closer to $x_i$ than $\frac\ve{16}$. 
To define elements $r^a_i$, 
}%
let $P(x_1,\dots,x_n)$ be the system of all equations of 
the form $x_i+x_j=x_k$ with unknowns $x_i,x_j,x_k$, where 
$1\le i,j,k\le n$ and in reality $y_i+y_j=y_k$. 
It follows from the choice of $D$ that this system has a 
solution $\ang{r_1,\dots,r_n}$ such that 
$r_i\in U^a_i\cap D$ for all $i$. 
In other words we have: $r_i+r_j=r_k$ whenever $y_i+y_j=y_k$. 
Let $r^a_i=r_i$. 
This ends the definition of $a\in A$ as required. 
(An extra care to guarantee that $U^a=\bigcup_{1\le i\le n}U^a_i$ 
is a proper subset of $\dH$ is left to the reader.)
\epf

\subsubsection{The ideal}
\label{z}

Let $\cZ$ be the set of all sets $X\sq A$ such that 
there is a finite set $u\sq\dH$ satisfying the 
following: 
for any $a\in X$ we have $u\not\sq U^a$. 

The idea of this ideal goes back to Solecki~\cite{so}, 
where a certain ideal over the set 
$\Om$ of all clopen sets $U\sq2^\dN$ of measure $\frac12$ 
(also a countable set) is considered. 
In our case the index set $A$ is somewhat more complicated. 


\ble
\label p
$\cZ$ is an ideal containing all finite sets\/ $X\sq A$, 
but\/ $A\nin\cZ$.
\ele
\bpf 
If $a\in A$ then the singleton $\ans a$ belongs to $\cZ$. 
Indeed by definition $U^a$ is a non-empty subset of $\dH$. 
Therefore there is a point $x\in \dH\bez U^a$. 
Then $u=\ans x$ witnesses $A\in \cZ$.  
To see that $\cZ$ is closed under finite unions, suppose 
that finite sets $u,v\sq\dH$ witness that resp.\ $X,Y$ 
belong to $\cZ$. 
Then $w=u\cup v$ obviously witnesses that $Z=X\cup Y\in\cZ$. 
Finally by Lemma~\ref{l*} for any finite 
$u=\ans{x_1,...,x_n}\sq\dH$ there is an element $a\in A$ 
such that $u\sq U^a$.
This implies that $A$ itself does not belong to $\cZ$.
\epf


\bpro
\label{p1}
$\cZ$ is an analytic ideal. 
If\/ $\dH$ is\/ \dd\sg compact then\/ $\cZ$ is\/ $\fs$.
\epro
\bpf
We claim that 
$X\in\cZ$ iff there are a natural $n$ and a partition 
$X=\bigcup_{1\le k\le n}X_k$ such that for any $k$ the set 
$X_k\sq A$ satisfies $\bigcup_{a\in X_k}U^a\ne\dH$. 
Indeed suppose that $X\in\cZ$ and this is witnessed by a 
finite set $u=\ans{x_1,\dots,x_n}\sq \dH$, that is,   
$u\not\sq U^a$ for all $a\in X$. 
It follows that $X=\bigcup_{1\le k\le n}X_k$, where 
$X_k=\ens{a\in X}{x_k\nin U^a}$. 
Clearly $x_k\nin \bigcup_{a\in X_k}U^a$. 
To prove the converse suppose that 
$X=\bigcup_{1\le k\le n}X_k\sq A$ and 
$\bigcup_{a\in X_k}U^a\ne\dH$ for all $k$. 
Let us pick arbitrary points 
$x_k\in\dH\bez\bigcup_{a\in X_k}U^a$ for all $k$.
Then $u=\ans{x_1,\dots,x_n}$ witnesses $X\in\cZ$, as required.

It easily follows that $\cZ$ is analytic. 

Now suppose that $\dH=\bigcup_{\ell\in\dN} H_\ell$, 
where all sets $H_\ell$ are compact.
Then the inequality $\bigcup_{a\in X_k}U^a\ne\dH$ is equivalent 
to $\sus\ell\:(H_\ell\not\sq\bigcup_{a\in X_k}U^a)$. 
And by the compactness, the non-inclusion 
$H_\ell\not\sq\bigcup_{a\in X_k}U^a$ is equivalent to 
the following statement: $H_\ell\not\sq\bigcup_{a\in X'}U^a$   
for every finite $X'\sq X_k$. 
Fix an enumeration $A=\sis{a_n}{n\in\dN}$. 
Put $A\res m=\ens{a_j}{j<m}$. 
Using K\"onig's lemma, we conclude that 
$X\in\cZ$ iff there exist natural $\ell,n$ such that 
for any $m$ there exists a partition 
$X\cap({A\res m})=\bigcup_{k<n}X_k$, 
where for every $k$ we have $H_\ell\not\sq\bigcup_{a\in X_k}U^a$. 
And this is a $\fs$ definition for $\cZ$.
\epf

\subsubsection
{The main result}

Here we prove Theorem~\ref{m}.
Define a Borel map $f:\dH\to\dH^A$ as follows. 
Suppose that $x\in\dH$ and $a\in A$, $n^a=n$.
If\/ $x\in U^a_i$, $1\le i\le n$, then put $f_a(x)=x\omi r^a_i$. 
($\omi$ in the sense of the group $\dH$.) 
If\/ $x\nin U^a$ then put simply $f_a(x)=\bbo$.   

Finally define $f(x)=\sis{f_a(x)}{a\in A}$.
Clearly $f$ is a Borel map.

The maps $f_a$ do not look like \hm s $\dH\to\dH$. 
Nevertheless their combination surprisingly turns out to be an 
approximate \hm!

\ble
\label{l1}
$f:\dH\to\dH^A$ is a Borel and\/ \ddz approximate \hm.
\ele
\bpf
Let $x,\,y\in\dH$ and $z=x\opl y$. 
Prove that the set 
$$
C_{xy}=\ans{a:f_a(x)\opl f_a(y)\ne f_a(z)}
$$ 
belongs to $\cZ$. 
We assert that this is witnessed by the set $u=\ans{x,y,z}$, 
that is, if $a\in C_{xy}$ then at least one of the 
points $x,y,z$ is not a point in $U^a$. 
Or, equivalently, if $a\in A$ and $x,y,z$ belong to $U^a$  
then $f_a(x)\opl f_a(y)=f_a(z)$. 

To prove this fact suppose that $a\in A$ and $x,y,z\in U^a$.
By definition, $U^a=U^a_1\cup\dots\cup U^a_n$, where $n=n^a$ 
and $U^a_i$ are disjoint 
rational balls in $\dH$. 
We have $x\in U^a_i$, $y\in U^a_j$, $z\in U^a_k$, where 
$1\le i,j,k\le n$. 
Then by definition 
$$
f_a(x)=x\omi r^a_i\,,\quad f_a(y)=y\omi r^a_j\,,\quad 
f_a(z)=z\omi r^a_k\,.
$$
Therefore $f_a(x)\opl f_a(y)=x\opl y\omi(s_i\opl s_j)$. 
(Here we clearly use the assumption that the group is 
abelian.)
We assert that $r^a_i\opl r^a_j=r^a_k$ --- then obviously 
$f_a(x)\opl f_a(y)=f_a(z)$ by the above, and we are done. 

Note that $z=x\opl y\in U^a$, hence 
$(U^a_i\opl U^a_j)\cap U^a\ne\pu$. 
We conclude that \ref{32} of Definition \ref{3a} fails. 
Therefore \ref{31} holds,   
$r^a_i\opl r^a_j=r^a_{k'}$ for 
some $k'$ and $(U^a_i\opl U^a_j)\cap U^a\sq U^a_{k'}$. 
But the set $(U^a_i\opl U^a_j)\cap U^a$ obviously contains $z$, and 
$z\in U^a_k$. 
It follows that $k'=k$, $r^a_{k'}=r^a_k$, $r^a_i\opl r^a_j=r^a_k$,  
as required.
\epF{Lemma}

\ble
\label{l2}
The approximate \hm\/
$f$ is not\/ \ddz approximable by a continuous \hm\/ 
$g:\dH\to\dH^A$.
\ele
\bpf
Assume towards the contrary  that $g:\dH\to\dH^A$ is 
a continuous \hm\ which \ddz approximates $f$. 
Thus if $x\in\dH$ then the set 
$\Da_x=\ans{a:f_a(x)\ne g_a(x)}$ belongs to $\cZ,$ 
where, as usual, $g_a(x)=g(x)(a)$. 
Note that all of these projection maps $g_a:\dH\to\dH$ are 
continuous group \hm s since such is $g$ itself.

Thus if $x\in\dH$ then $\Da_x\in\dZ$, and hence there is a 
finite set $u_x\sq D$ satisfying the following: 
if $a\in A$ and $u_x\sq U^a$ then $a\nin \Da_x$, that is, 
$f_a(x)=g_a(x)$.
Put
\dm
X_u=\ens{x\in\dH}
{\kaz a\in A\:(u\sq U^a\imp f_a(x)=g_a(x))}
\dm
for every finite $u\sq D$. 
These sets are Borel since so are maps $f,g$ 
(and $g$ even continuous). 
Moreover $\dH=\bigcup_{u\sq D\;\text{finite}}X_u$ since every 
$x\in\dH$ belongs to $X_{u_x}$.
Thus at least one of the sets $X_u$ is not meager, therefore, 
is comeager on a certain rational ball $B\sq\dH$.
Fix $u$ and $B$. 
By definition for comeager-many $x\in B$ and all  
$a\in A$ satisfying $u\sq U^a$ we have $f_a(x)=g_a(x)$.

Arguing as in the proof of Lemma~\ref{l*}, we obtain an 
element $a\in A$ satisfying the following properties: 
$u\sq U^a$, $U^a\cap B\ne\pu$, but the set $B\bez U^a$ 
is non-empty and moreover is not dense in $B$. 
Fix such $a$.
Thus there exists a non-empty rational ball $B'\sq B$ that 
does not intersect $U^a$. 
By definition $f_a(x)=\bbo$ for all $x\in B'$, and hence 
$g_a(x)=\bbo$ for comeager-many $x\in B'$ by the choice of $B$. 
We conclude that $g_a(x)=\bbo$ for all $x\in B$ in general, 
because $g$ is continuous.         

Now, let $n^a=n$. 
Then $U^a=U^a_1\cup\dots\cup U^a_n$. 
Recall that the intersection $B\cap U^a$ of two open sets is 
non-empty by the choice of $a$. 
It follows that there exists an index $i$, $1\le i\le n$, and 
a non-empty rational ball $B''\sq B\cap U^a_i$.
Then by definition $f_a(x)=x\omi r$ for all $x\in B''$, 
where $r=r^a_i$. 
Therefore $g_a(x)= x\omi r$ for comeager-many $x\in B''$, 
and then $g_a(x)= x\omi r$ for all $x\in B''$ since $g$ is 
continuous.

To conclude, $g_a$, a continuous group \hm, is constant $\bbo$ 
on a non-empty open set $B'$, and is bijective on another 
non-empty open set $B''$. 
But this cannot be the case.
\epF{Lemma}
 
Lemmas \ref{l1} and \ref{l2} complete the proof of 
Theorem~\ref{m}.

\end{document}